\documentclass[11pt]{article}
\usepackage{}
\usepackage{amsfonts}
\usepackage{amssymb}
\usepackage{amsmath,amssymb,graphicx}
\usepackage{epsfig}
\newtheorem{theorem}{Theorem}[section]
\newtheorem{lemma}[theorem]{Lemma}
\newtheorem{proposition}[theorem]{Proposition}

\newcommand{\diver}{\mathrm{div}\,}
\newcommand{\ma}{\iota}
\newcommand{\mj}{\textbf{j}}
\newcommand{\mcs}{\mathcal{S}}
\newcommand{\kz}{\kappa(z)}
\newcommand{\nn}{\nonumber}

\newcommand{\be}{\begin{eqnarray}}
\newcommand{\ee}{\end{eqnarray}}
\newcommand{\ba}{\begin{array}}
\newcommand{\ea}{\end{array}}
\newcommand{\ben}{\begin{eqnarray*}}
\newcommand{\een}{\end{eqnarray*}}

\newcommand{\pa}{\partial}

\begin{document}

\hoffset -34pt

\title{{\large \bf Variants of Ahlfors' lemma and properties of the logarithmic potentials
}
\thanks{\mbox{Keywords. Gaussian curvatures, SK-metrics, logarithmic potentials.}}}
\author{\normalsize  Tanran Zhang }
\date{}
\maketitle \baselineskip 21pt
\noindent

\begin{minipage}{138mm}
\renewcommand{\baselinestretch}{1} \normalsize
\begin{abstract}
{As a special class of conformal metrics with negative curvatures, SK-metrics play a crucial role in metric spaces. This paper concerns the variants of Ahlfors' lemma in an SK-metric space and gives a higher order derivative formula for the logarithmic potential function, which can be applied for the estimates near the singularity of a conformal metric with negative curvatures. }
\end{abstract}
\end{minipage}\\
\\
\renewcommand{\baselinestretch}{1} \normalsize
\section{Introduction}

Let $\mathbb{C}$ be the complex plane, $\mathbb{D}$ be the open unit disk in $\mathbb{C}$ and $\mathbb{D^*}:=\mathbb{D} \setminus \{0\}$ be the punctured unit disk. Let $\mathcal{S}$ denote a Riemann surface. Let $A$ be the basic family consisting of homeomorphisms $\sigma$ defined on the plane domains of $\mathbb{C}$ into $\mathcal{S}$ which defines the conformal structure of $\mathcal{S}$. In our discussion we take the linear notation for a conformal metric $ds=\lambda(z)|dz|$ and set the density function  $\lambda(z)$ to be positive on $\mathcal{S}$. However, in Heins' paper \cite{Heins}, he let the function $\lambda(z)$ be non-negative for a conformal metric $\lambda(z)|dz|$ on $\mathcal{S}$, and then defined the SK-metrics. \\
\par Let $\mathbb{P}$ denote the Riemann sphere $\mathbb{C}\cup \{\infty\}$ with its canonical complex structure and a subdomain $\widetilde{\Omega}\subset\mathbb{P}$. For a point $p \in \widetilde{\Omega}$, let $z$ be the local coordinates such that $z(p)=0$. We say that a conformal metric $\lambda(z)|dz|$ on the punctured domain $\Omega:=\widetilde{\Omega} \backslash \{p\}$ has a \emph{singularity of order $\alpha\leq 1$ at the point $p$}, if, in local coordinates $z$,
\be \label{singularity}
\log\lambda(z)=\left\{
\begin{array}{ll}
-\alpha\log|z|+\textit {O}(1) & \mbox{if\ }\ \alpha < 1 \\
-\log|z|-\log\log(1/|z|)+\textit {O}(1)&\mbox{if\ }\ \alpha=1 \end{array} \right. \ee
as $z\rightarrow 0$, with $\textit {O}$ being the Landau symbol throughout our study. Let $M_u(r):=\sup_{|z|=r}u(z)$ for a real-valued function $u(z)$ defined in a punctured neighborhood of $z=0$ and call
\be\label{order of u}
\alpha(u):=\lim_{r\rightarrow 0^+}\frac{M_u(r)}{\log(1/r)}
\ee
the \emph{order of $u(z)$} if this limit exists. For $u(z):=\log \lambda(z)$, $\alpha(u)$ in (\ref{singularity}) and \eqref{order of u} are equivalent. We call the point $p$ a \emph{conical singularity} or \emph{corner of order} $\alpha$ if $\alpha< 1$ and a \emph{cusp} if $\alpha=1$. The \emph{generalized Gaussian curvature $\kappa_{\lambda}(z)$} of the density function $\lambda(z)$ is defined by
\be \label{general kappa}
\kappa_{\lambda}(z)=-\frac{1}{\lambda(z)^2}{\liminf_{r\rightarrow 0}\frac{4}{r^2}(\frac 1 {2\pi}\int_0^{2\pi}\log\lambda(z+re^{it})dt-\log\lambda(z))}.
\ee
We say a conformal metric $\lambda(z)|dz|$ on a domain $\Omega\subseteq\mathbb{C}$ is \emph{regular}, if its density $\lambda(z)$ is of class $C^2$ in $\Omega$. We will show that, if $\lambda(z)|dz|$ is a regular conformal metric, then (\ref{general kappa}) is equivalent to
\be \label{special kappa}
\kappa_{\lambda}(z)=-\frac{\Delta\log\lambda(z)}{\lambda(z)^2},
\ee
where $\Delta$ denotes the Laplace operator. We will discuss details in Lemma \ref{mean value}. It is well known that, if $a<\kappa_{\lambda}(z)<b<0$, the metric $\lambda(z)|dz|$ only has corners or cusps at isolated singularities, see \cite{McOwen1}. If $\kappa_{\lambda}(z)\geq 0$ and the energy is finite, then only corners occur, see \cite{Rothequan,Troyanov1}, also \cite{Rothhyper}.
\par We have a fact that the Gaussian curvature is a conformal invariant. Let $\lambda(z)|dz|$ be a conformal metric on a plane domain $D$ and $f: W\rightarrow D$ be a holomorphic mapping of a Riemann surface $W$ into $D$. Then the metric $ds=f^{*}\lambda(w)|dw|$ on $W$ induced by $f$ from the original metric $\lambda(z)|dz|$ is called the \emph{pullback of $\lambda(z)|dz|$} and defined by
\be \label{pullback}
ds=f^{*}\lambda(w)|dw|:=\lambda(f(w))|f'(w)||dw|.
\ee
It is evident that $f^{*}\lambda(w)|dw|$ is a conformal metric on $W \backslash\{\mbox{critical points of}\ f\}$ with Gaussian curvature
\be
\kappa_{f^{*}{\lambda}}(w)=\kappa_{\lambda}(f(w)). \nn
\ee
Using this conformal invariance, we can only consider one Riemann surface with the conformal metric all over its conformal equivalence class.
\par The paper is divided into two parts. Section 2 is devoted to a class of conformal metrics, SK-metrics. We begin with the definition and then discuss the maximum principle for these metrics. Section 3 is about the potential theory, which is a main tool we use to study the local term of a conformal metric near the singularity.

\section {SK-metrics}
\setcounter{equation}{0}

For a topological space $X$, a function $u: X \rightarrow [-\infty, \infty)$ is called \emph{upper semi-continuous} if the set $\{z \in X: u(z)<\alpha\}$ is open in $X$ for each real number $\alpha$. If $\varphi$ is a \emph{uniformizer} of $\mcs$, i.e. $\varphi$ is a univalent conformal map of a plane domain into $\mcs$, then the conformal density function $\lambda$ on $\mcs$ can be extended to $\mcs\cup \{\varphi\}$ such that the extension is a conformal metric relative to $\mcs\cup \{\varphi\}$. We denote the image of $\varphi$ with respect to this extension by $\lambda_{\varphi}$ and call it the \emph{$\varphi$-scale of $\lambda$}. For two uniformizers $\varphi$ and $\psi$ of $\mcs$, $\varphi(a)=p,$\ $\psi(b)=p,$\ $a,b \in \mcs,$\ $p \in \mathbb{C},$\ if $\lambda_{\varphi}$ and $\lambda_{\psi}$ are upper semi-continuous at $a$ and $b$ respectively, we say $\lambda$ is upper semi-continuous at $p$. If $\lambda$ is upper semi-continuous at every point of $\mcs$, we say $\lambda$ is upper semi-continuous on $\mcs$.

Now we give the definition for SK-metrics. The concept of SK-metrics was given by Heins in \cite{Heins}, but its initial idea came from Ahlfors in \cite{Ahlforslemma} and \cite{Ahlforsbook}. Heins defined the SK-metrics by means of the Gaussian curvature and proved a more general maximum principle as a variant of Ahlfors' lemma. According to Heins, the terminology of "SK-metric" is partly because its (Gaussian) "curvature is subordinate to $-4$", see \cite[p.3]{Heins}. We call an upper semi-continuous metric $\lambda(z)|dz|$ on $\mathcal{S}$ an \emph{SK-metric} if its Gaussian curvature is bounded above by $-4$ at those points $z$ in $\mathcal{S}$ where $z$ satisfies $\lambda(z)>0$. A complete metric with the negative constant Gaussian curvature is called the hyperbolic metric. The hyperbolic metric on the unit disk $\mathbb{D}$ is defined by
\be \label{hyperbolic metric}
d\sigma=\lambda_{\mathbb{D}}(z)|dz|=\frac {|dz|}{1-{|z|}^2}
\ee
with the constant Gaussian curvature $-4$ and it is an extremal SK-metric. For SK-metrics, there is a generalization of the maximum principle mentioned by Ahlfors [2, Theorem A] and Heins [4, Theorem 2.1], which claims that the hyperbolic metric is the unique maximal SK-metric on $\mathbb{D}$.
\begin{theorem} $\mathrm{[2]\ (Ahlfors'\ lemma)}$ \label{maximum of Ahl}
\textsl{Let $d\sigma$ be the hyperbolic metric on $\mathbb{D}$ given in (\ref{hyperbolic metric}) and $ds$ be the metrics on $\mathbb{D}$ induce by an SK-metric on some Riemann surface $W$. If the function $f(z)$ is analytic in $\mathbb{D}$, then the inequality
\be
ds \leq d\sigma \nn
\ee
will hold throughout the circle. }
\end{theorem}
\par The following result stated on the corresponding Riemann surface. It is a variant of Theorem \ref{maximum of Ahl}.
\begin{theorem} $\mathrm{[4]}$ \label{maximum of Heins}
\textsl{Suppose that $W$ is a relatively compact domain of $\mcs$ and that $\lambda(w)$ is an SK-metric on $W$, $\mu(w)$ is a pullback on $W$ of $\lambda_{\mathbb{D}}(z)$ defined in (\ref{hyperbolic metric}). If for all $p \in \pa W$,
$$ \displaystyle \limsup_{w \rightarrow p} \frac{\lambda(w)}{\mu(w)} \leq 1, $$
then throughout the boundary $\pa W$, it holds that $$\lambda(w) \leq \mu(w). $$}
\end{theorem}

Heins used condition (\ref{general kappa}) to define SK-metrics and he mentioned the equivalence between (\ref{general kappa}) and (\ref{special kappa}) for SK-metrics in one word, see \cite[(1.4)]{Heins}. Here we present it in detail.
\begin{lemma} \label{mean value}
\textsl{Suppose $\Omega \in \mathbb{C}$ is a domain. If a function $u$ is of class $C^2(\Omega)$, then for $z \in \Omega$, we have
$$\lim_{r\rightarrow 0}\frac{4}{r^{2}} \left(\frac{1}{2\pi}\int^{2\pi}_{0}u(z+re^{it})dt-u(z)\right)=
\Delta u(z).$$}
\end{lemma}
\textbf{Proof.} Since $u$ is of  $C^2$, using Taylor's expansion of $u(z)$ at $z_0 \in \Omega$, \\
$$u(z_0+z)=u(z_0)+u_x(z_0) \cdot x+u_y(z_0) \cdot y + \frac{1}{2}u_{xx}(z_0) \cdot {x^{2}}+ \frac{1}{2}u_{yy}(z_0) \cdot {y^{2}} +u_{xy}(z_0) \cdot xy + \varepsilon (z),$$
where $\varepsilon (z)\rightarrow 0$ as $z \rightarrow 0$ and $z=x+yi$,
we have
$$\frac{1}{2\pi}\int^{2\pi}_{0}u(z_0+re^{it})dt-u(z_0)=
\frac{r^2}{4}\left(u_{xx}(z_0)+u_{yy}(z_0)\right)
+ \frac{1}{2\pi}\int^{2\pi}_{0}\varepsilon (z)dt.$$
As $r=|z|\rightarrow 0$,
$$\frac{1}{r^2}\int^{2\pi}_{0}\varepsilon (z)dt=
\int^{2\pi}_{0}\frac{\varepsilon (z)}{r^2}dt=0,$$ then
$$\lim_{r\rightarrow 0}\frac{4}{r^{2}}\left(\frac{1}{2\pi}\int^{2\pi}_{0}u(z+re^{it})dt-u(z)\right)=
u_{xx}(z_0)+u_{yy}(z_0)
=\Delta u(z_0)$$
as required. \hfill $\Box$
\par The following theorem offers us a simple way to construct a new SK-metric on a plane domain by means of the maximum principle. This is called the "gluing lemma".
\begin{theorem} $\mathrm{[7]}$
\textsl{Let $\lambda(z)|dz|$ be an SK-metric on a domain $G \in \mathbb{C}$ and let $\mu(z)|dz|$ be an SK-metric on a subdomain $U$ of $G$ such that the "gluing condition"
$$\limsup_{U \ni z\rightarrow \xi}\mu(z) \leq \lambda(\xi)$$
holds for all $\xi \in \pa U \cap G$. Then $\sigma(z)|dz|$ defined by
\ben
\sigma(z):=\left\{
\begin{array}{ll}
\max\{\lambda(z), \mu(z)\} & \mbox{for\ }\ z \in U, \\
\lambda(z) &\mbox{for\ }\ z \in G \setminus U
\end{array} \right.
\een
is an SK-metric on $G$.}
\end{theorem}
\par We end this section with the discussion on SK-metrics in the punctured domain. On the punctured unit disk $\mathbb{D}^*$, the hyperbolic metric is expressed by
\be \label{poin metric punc}
\lambda_{\mathbb{D}^*}(z)|dz|=\frac{|dz|}{2|z|\log(1/|z|)}
\ee
with the constant curvature $-4$. From the definition (\ref{singularity}) of the singularity and its order, we know that the metric (\ref{poin metric punc}) has order $1$ at the origin. In any punctured disk, Kraus, Roth and Sugawa gave the expression of the hyperbolic metric which has a singularity at the origin of order $\alpha <1$ in \cite{Rothhyper} without any detailed discussion. Now we give a complete presentation of the proof.

\begin{theorem} $\mathrm{[7]}$ \label{maximal}
\textsl{ For $R>0$, let $D_R^*:=\{z: 0<|z|<R \}$ and
\ben
\lambda_{\alpha,R}(z):=\left\{
\begin{array}{ll}
\displaystyle \frac{1-\alpha}{2|z|\sinh \left((1-\alpha)\log {\frac{R}{|z|}}\right)} & \mbox{if\ }\ \alpha<1, \\
\displaystyle \frac{1}{2|z|\log \frac{R}{|z|}} &\mbox{if\ }\ \alpha=1
\end{array} \right.
\een
for $z \in D_R^*$, then for an arbitrary SK-metric $\sigma(z)$ on $D_R^*$ which has a singularity at $z=0$ of order $\alpha$, we have $\sigma(z)\leq\lambda_{\alpha,R}(z)$.}
\end{theorem}
\textbf{Proof.} We consider the case $\alpha<1$. First, choose an arbitrary $0 < R_0 < R$, consider $\lambda_{\beta,R_0}(z)$ on $0<|z|<R_0$ for $\alpha<\beta<1$, and let $u(z):=\log\sigma(z)$, $v(z):=\log\lambda_{\beta,R_0}(z)$, $E:=\{z:0<|z|<R_0, u(z)>v(z)\}$.
\par Now we have the assertion that $0 \notin \overline{E}$.
Since $\sigma(z)|dz|$ and $\lambda_{\beta,R_0}(z)|dz|$ both have a singularity at $z=0$ with order $\alpha$, $\beta$
respectively, then $$v(z)=-\beta\log|z|+O(1), \, \,  u(z)=-\alpha\log|z|+O(1),$$ so
$u-v=(\beta-\alpha)\log|z|+O(1)$. Since $u-v\rightarrow -\infty$ as $z\rightarrow 0$, then on a sufficiently small neighborhood of $z=0$, $u-v<0$ holds, thus
$0 \notin \overline{E}$.
\par Similarly, we have $\partial E \cap \{z: 0<|z|<R_0 \}=\emptyset$, because $v \rightarrow + \infty$ as $|z| \rightarrow R_0$, and $u$ is bounded in $\{z: 0<|z|<R_0\}$.
\par Then consider the curve $|z|=R_0$. It is clear that $v(z)$ and $u(z)$ satisfy
$$ \displaystyle \lim_{r\rightarrow 0}\frac{4}{r^{2}}\left(\frac{1}{2\pi}\int^{2\pi}_{0}v(z+re^{it})dt-v(z)\right)=e^{2v},$$
and
$$ \displaystyle \liminf_{r\rightarrow 0}\frac{4}{r^{2}}\left(\frac{1}{2\pi}\int^{2\pi}_{0}u(z+re^{it})dt-u(z)\right)\geq e^{2u}$$
by Lemma \ref{mean value}.
Thus $$ \displaystyle \liminf_{r\rightarrow 0}\frac{4}{r^{2}}\left(\frac{1}{2\pi}\int^{2\pi}_{0}\left(u(z+re^{it})-v(z+re^{it})\right)dt-
\left(u(z)-v(z)\right)\right)\geq e^{2u}-e^{2v},$$
which is positive on $E$. By definition of limit inferior, we have for $z \in E$
$$\frac{1}{2\pi}\int^{2\pi}_{0}\left( u(z+re^{it})-v(z+re^{it})\right) dt-\left( u(z)-v(z)\right) > 0,$$
therefore, $$u(z)-v(z)\leq \frac{1}{2\pi}\int^{2\pi}_{0}\left(u(z+re^{it})-v(z+re^{it})\right)dt.$$
Now we recall the definition of subharmonic functions. Let $\Omega$ be an open
subset of $\mathbb{C}$. A function $u: \Omega \rightarrow [-\infty,
\infty)$ is called subharmonic if $u$ is upper semi-continuous and
satisfies the local sub-mean inequality, i.e. given $z \in \Omega$, there
exists $\rho>0$ such that
\be \label{circumferential mean}
u(z) \leq \frac{1}{2\pi}
\int_{0}^{2\pi}u(z+re^{it})dt
\ee
for $0 \leq r < \rho$.
If we adopt the definition as above, then $u-v$ is subharmonic on $E$, hence $u-v$ has no maximum in $E$ and $u-v$
approaches its least upper bound on a sequence tending to $\partial E$. A contradiction. So $E=\emptyset$.
\par Finally, letting $R_0 \rightarrow R$ and $\beta \rightarrow \alpha$ gives the maximality of $\lambda_{\alpha,R}(z)$ for $\alpha<1$. According to Kraus, Roth and Sugawa, for the case $\alpha=1$ this expression has to be interpreted in the limit sense $\alpha\nearrow 1$ to obtain $\lambda_{1,R}(z)$, i.e.
$$\lambda_{1,R}(z)=\lim_{\alpha\nearrow 1}\lambda_{\alpha,R}(z)= \frac{1}{2|z|\log \frac{R}{|z|}}$$ \hfill $\Box$\\[3mm]
\textbf{Remark.} The righthand side of (\ref{circumferential mean}) is called the circumferential mean of function $u$. Heins used it to describe the curvature in the definition of SK-metrics in \cite{Heins} with $\rho=1$ and $z=0$.\\

\section{Potential theory}
\setcounter{equation}{0}

Generally speaking, the SK-metric is defined by the fact that its Gaussian curvature no greater than some negative constant. So the maximum principle for the SK-metric is common and useful. After a combination with PDEs, the asymptotic behavior of a metric has something in-between with the local properties of the solution to the corresponding PDE. We can consider the curvature equation
\be\label{general equation}
\Delta u=-\kz e^{2u},
\ee
where $\kz$ is known, then the definition of SK-metrics is fruitful in the case that the curvature function $\kz$ is strictly negative and H\"{o}lder continuous in $\mathbb{D}^*$, see \cite{Rothbehaviour}, also \cite{Nitsche} for details. For an SK-metric $\lambda(z)$ on $\mathbb{D}^*$, regarding $u:=\log \lambda$ as a solution to the equation (\ref{general equation}), the global properties of $u$ have been well known by means of the study on equation (\ref{general equation}). However, near the singularity $z=0$, the local properties are still not explicit. We can employ a way related to partial differential equations to investigate the asymptotic behavior of $u$ near the origin. Potential theory is a powerful tool in our case. In this section, we present a formula of the higher order derivatives for the logarithmic potential, and give an asymptotic estimate for $u$ near the origin without any proof as an application of potential theory. We only refer to the logarithmic potential.
\par We identify $\mathbb{C}$ with $\mathbb{R}^2$, and write $z=x_1+ix_2$, $\zeta=y_1+iy_2$. Set $0<r\leq 1$ and denote $D_R:=\{z \in {\mathbb{C}}:|z|<R\}$, ${D_R}^*:=D_R\backslash \{0\}$ for a positive number $R$.
For a bounded, integrable function $f(z)$ defined on a domain $\Omega \subseteq \mathbb{C}$, the integral $$\frac {1}{2\pi}\iint _{\Omega}\log|z-\zeta| f(\zeta) d\sigma_{\zeta}$$ is called the \emph{logarithmic potential of $f$}, where $d\sigma_{\zeta}$ is the area element.
The H\"{o}lder spaces $C^{n, \nu}(D_R)$ are defined as the subspaces of $C^n(D_R)$ consisting of functions whose $n-$th order partial derivatives are locally H\"{o}lder continuous with exponent $\nu$ in $D_R$, $0<\nu \leq 1$. Then the following proposition for the first and the second order derivatives of the logarithmic potential is valid.

\begin{proposition} $\mathrm{[3,\ 6]}$ \label{newton potential}
\textsl{ Let $f: D_r\rightarrow\mathbb{R}$ be a locally bounded, integrable function in $D_r$ and $\omega$ be the logarithmic potential of $f$. Then $\omega \in C^1({D_r})$ and for any $z=x_1+ix_2 \in D_r$,
$$\frac{\partial}{\partial x_j}\omega(z)=\frac 1 {2\pi}\iint_{D_r}\frac {\partial}{\partial x_j}\log|z-\zeta| f(\zeta) d\sigma_{\zeta}$$ for $j\ \in \{1,2\}$.\\
If, in addition, $f$ is locally H\"{o}lder continuous with exponent $\nu \leq 1$, then $\omega\in C^2(D_r)$,\ $\Delta \omega=f$ in $D_r$ and for $z \in D_r$,
\ben
\frac{\partial^2}{\partial x_l\partial x_j}\omega(z)&=&\frac 1 {2\pi}\iint_{D_R}\frac{\partial^2}{\partial x_l\partial x_j}\log|z-\zeta|\left(f(\zeta)-f(z)\right) d\sigma_{\zeta} \\
&&-\frac{1}{2\pi}f(z)\int_{\partial D_R}\frac{\partial}{\partial x_j}\log|z-\zeta|N_l(\zeta)|d\zeta|,
\een
where $N(\zeta)=(N_1(\zeta),N_2(\zeta))$ is the unit outward normal at the point $\zeta \in\partial D_R$,\ $R>r$ such that the divergence theorem holds on $D_R$ and $f$ is extended to vanish outside of $D_r$.}
\end{proposition}

\par There is a similar proposition for higher order derivatives of the logarithmic potential. Define a multi-index $\mj=(j_1, j_2)$,  $|\mj|=j_1+j_2$, $j_1, j_2=0,1,2,\ldots \,$. For $z=x_1+ix_2$, denote
$$ \displaystyle \frac{\pa}{\pa x_1}=\pa_1,\ \frac{\pa}{\pa x_2}=\pa_2,\ \partial^{\mj}=\pa^{j_1}_1 \pa^{j_2}_2\ \mathrm{and}\ \frac{\pa^j}{\pa \zeta}=\frac{\pa^{j_1}}{\pa y_1}\frac{\pa^{j_2}}{\pa y_2}.$$
Let $e_{\tau}=(0,1)$ or $(1,0)$ for $\tau=1,2,\ldots \,$. Then $\mj$ can be expressed in the form $e_1+e_2+\cdots+e_n$. We define two vectors $\theta_{\tau}:=e_1+ \cdots +e_{\tau}$, $\phi_{\tau}:=e_{\tau+2}+ \cdots +e_n$ for $\tau=1, \ldots, \,n-1$ where $\phi_{n-1}=(0,0)$, so $\mj$ has a decomposition as $\mj=\theta_{\tau}+e_{\tau+1}+\phi_{\tau}$. Write $\zeta=y_1+iy_2$ and set
\begin{eqnarray*}
          P_n[f](z,\zeta):=\left\{
            \begin{array}{ll}
            \vspace*{2mm}
             \displaystyle \sum_{|\ma|\leq n} \frac{(\zeta-z)^{\ma}{\pa}^{\ma}}{\ma !}f(z) & \mbox{if\ }\ n \geq 1 \\
             f(z) & \mbox{if\ }\ n=0,
            \end{array} \right.
\end{eqnarray*}
where $\ma$ is a multi-index, $\ma=(\iota_1, \iota_2)$, $(\zeta-z)^{\ma}=(y_1-x_1)^{\iota_1}(y_2-x_2)^{\iota_2}$, $\ma !=\iota_1!\iota_2!$. We have the following recurrent formula for $P_n[f](z,\zeta)$.
\begin{lemma} \label{recurrent formula}
\textsl{For $P_n[f](z,\zeta)$ defined as above, then
$$\frac{\pa^{e}}{\pa \zeta}P_n[f](z,\zeta)=P_{n-1}[\pa^e f](z,\zeta)$$
holds for $e=(0,1)$ or $(1,0)$. }
\end{lemma}
\textbf{Proof.} We take the case $e=(1,0)$ as an example, when $e=(0,1)$ it is similar. Let $n \geq 1$. Then
\begin{eqnarray*}
&&\frac{\pa}{\pa y_1}P_n[f](z,\zeta)  \\
&=&\frac{\pa}{\pa y_1} \sum_{{\ma_1+\ma_2 \leq n} \atop {0 \leq \ma_1 \leq n}} \frac{(y_1-x_1)^{\ma_1}(y_2-x_2)^{\ma_2}}{\ma_1 ! \ma_2 !}{\pa_1}^{\ma_1}{\pa_2}^{\ma_2}f(z) \\
&=& \sum_{{\ma_1+\ma_2 \leq n} \atop {1 \leq \ma_1 \leq n}} \frac{(y_1-x_1)^{\ma_1-1}(y_2-x_2)^{\ma_2}}{(\ma_1-1) ! \ma_2 !}{\pa_1}^{\ma_1}{\pa_2}^{\ma_2}f(z) \\
&=& \sum_{{(\ma_1-1)+\ma_2 \leq n-1} \atop {0 \leq \ma_1-1 \leq n}} \frac{(y_1-x_1)^{\ma_1-1}(y_2-x_2)^{\ma_2}}{(\ma_1-1) ! \ma_2 !}{\pa_1}^{\ma_1-1}{\pa_2}^{\ma_2}\, {\pa_1}f(z)\\
&=& \sum_{\ma_1+\ma_2 \leq n-1} \frac{(y_1-x_1)^{\ma_1}(y_2-x_2)^{\ma_2}{\pa_1}^{\ma_1}{\pa_2}^{\ma_2}}{\ma_1 ! \ma_2 !} {\pa_1}f(z)\\
&=& \sum_{|\ma|\leq n-1} \frac{(\zeta-z)^{\ma}{\pa}^{\ma}}{\ma !} \pa_1 f(z)\\
&=& P_{n-1}[\pa_1 f](z,\zeta).
\end{eqnarray*}
\hfill  $\Box$ \\
Using the multi-index notation, the Taylor expansion of $f$ can be written in short.
\begin{theorem}  \label{Taylor theorem} $\mathrm{[cf.\ 1]}$
\textsl{If $f(\zeta)$ is analytic in a domain $\Omega \in \mathbb{C}$, containing the point $z$, it is possible to write
$$f(\zeta)=\sum^{n}_{t=0}\frac{f^{(t)}(z)}{t!}(\zeta-z)+R_{n+1}(z,\zeta),$$
where $R_{n+1}(z,\zeta)$ is the error term and $R_{n+1}(z,\zeta)=f_{n+1}(z)(\zeta-z)^{n+1}$ with $f_{n+1}(z)$ analytic in $\Omega$. This expression is equivalent to 
\be \label{Taylor expansion}
f(\zeta)=P_n[z](z,\zeta)+R_{n+1}(z,\zeta),
\ee
with $R_{n+1}(z,\zeta)$ as above.}
\end{theorem}
\textbf{Remark.} If $f \in C^{n,\nu}(\Omega)$ with $0<\nu \leq 1$,\ $n\geq 1$, and the H\"{o}lder continuity is a local property, then the error term $R_{n+1}(z,\zeta)$ satisfies
\be \label{error term}
R_{n+1}(z, \zeta)=\textit{O}(|z-\zeta|^{\nu+n}).
\ee
\par On the basis of Lemma \ref{recurrent formula}, we can present the analogue of Proposition \ref{newton potential} as follows.
\begin{proposition} \label{new higher prop}
\textsl{ Let $r<1$, $f: D_r\rightarrow\mathbb{R}$, $f \in C^{n-2,  \nu}(D_r)$ with $0<\nu \leq 1$, $n\geq 3$ and $\omega$ be the logarithmic potential of $f$. Then
$\omega(z)\in C^{n}(D_r)$ and for a multi-index $j$,\ $|\mj|=n$,
\be \label{new higher}
\partial^{\mj}\omega(z)&=&\frac 1 {2\pi}\iint_{D_R}\partial^{\mj} \log|z-\zeta| \cdot \left(f(\zeta)-P_{n-2}[f](z,\zeta)\right)d\sigma_{\zeta} \nn \\
&& -\frac{1}{2\pi} \sum^{n-1}_{\tau=1} \int_{\partial D_R}\partial^{\theta_{\tau}} \log|z-\zeta| \cdot P_{\tau-1}[\partial^{\phi_{\tau}}f](z,\zeta)\cdot \langle N(\zeta),e_{\tau+1} \rangle |d\zeta|,
\ee
where $N(\zeta)=(N_1(\zeta),N_2(\zeta))$ is the unit outward normal at the point  $\zeta \in \partial D_R$, $\langle \ , \  \rangle$ is the inner product, $R>r$ such that the divergence theorem holds on $D_R$ and the function $f$ is extended to vanish outside of $D_r$.}
\end{proposition}
\par We need the following Divergence Theorem for the proof. For a point $z=(x_1, x_2)$, a vector field $\mathbf{w}(z)=(w_1(z),w_2(z))$ and a function $u(z)$, denote $$\diver \mathbf{w}=\frac{\pa w_1}{\pa x_1}+\frac{\pa w_2}{\pa x_2}=\textrm{divergence\ of}\ \mathbf{w},$$
$$D u=(\pa_1 u, \pa_2 u)=\textrm{gradient\ of}\ u,$$
then $\Delta u=\diver Du$.
\begin{theorem}  \label{divergence theorem}
$\mathrm{[cf.\ 3]\ (Divergence\ Theorem)\ }$ \textsl{Let $\Omega$ be a bounded domain with $C^1$ boundary $\partial \Omega$, for any vector field $\mathbf{w}$ in $C^0(\bar{\Omega})\cap C^1(\Omega)$, we have
\be \label{divergence equality}
\iint_{\Omega} \diver \mathbf{w} d\sigma_z=\int_{\partial\Omega} \langle N(z), \mathbf{w} \rangle \,d|z|,
\ee
where $ \langle \ , \  \rangle $ is the inner product. }
\end{theorem}
In (\ref{divergence equality}) we select $\mathbf{w}(z)=v(z)\,Du(z)$, then
\be \label{in between}
\iint_{\Omega} Du\, Dv \,d\sigma_z+ \iint_{\Omega} v\, Du\,d\sigma_z=\int_{\partial\Omega} v\, \langle Du, N(z) \rangle \,d|z|.
\ee
Since we only need one $\partial_{m}$ for $m=1,2$, we can fix the the other component $x_{3-m}$ in (\ref{in between}) and relabel $u$, we obtain the following \textit{Green's (first) identity}:
\be \label{Green's identity}
\iint_{\Omega} u\, \pa_m v \,d\sigma_z+ \iint_{\Omega} v\, \pa_m u \,d\sigma_z=\int_{\partial\Omega} u v\, N_m(z)\, d|z|.
\ee
\\[6mm]
\textbf{Proof of Proposition \ref{new higher prop}.} Let
\be \label{ujz}
u_j(z)&=&\frac 1 {2\pi}\iint_{D_R}\partial^{\mj} \log|z-\zeta| \cdot \left(f(\zeta)-P_{n-2}[f](z,\zeta)\right)d\sigma_{\zeta} \nn \\
&& -\frac{1}{2\pi} \sum^{n-1}_{\tau=1} \int_{\partial D_R}\partial^{\theta_{\tau}} \log|z-\zeta| \cdot P_{\tau-1}[\partial^{\phi_{\tau}}f](z,\zeta)\cdot \langle N(\zeta),e_{\tau+1} \rangle |d\zeta|.
\ee
Note that
\be \label{logn}
\left| \partial^{\mj}\log|z-\zeta| \right|
\leq \frac{n!}{|z-\zeta|^{n}},
\ee
for $n=|j|$, and $\log|z-\zeta|$ is harmonic for $\zeta \neq z$, then by the local H\"{o}lder continuity of $f$ in $D_r$, the function $u_j(z)$ is well defined.\\
Now we can employ induction. Since Proposition \ref{new higher prop} has been obtained already, and $j$ has the decomposition $j=\theta_{n-1}+e_n$,\ we may assume that the formula (\ref{new higher}) is true for $\theta_{n-1}$. Fix a function $\eta(t) \in C^{n-1}(\mathbb{R})$ such that $0 \leq \eta \leq 1$,\ $0 \leq \eta^{(n-1)} \leq 2$,\ $\eta(t)=0$ for $t\leq1$,\ $\eta(t)=1$ for $t\geq2$,\ and set $$\eta_{\varepsilon}:=\eta(\frac{|z-\zeta|}{\varepsilon}), \ L:=\frac{1}{2\pi}\log|z-\zeta|.$$ Note that $\eta_{\varepsilon}$ and $L$ are both skew symmetric with respect to $x_1$ and $y_1$,\ $x_2$ and $y_2$. Then
\be \label{skew symmety}
\pa^e L\eta_{\varepsilon}=-\frac{\pa^e}{\pa \zeta} L\eta_{\varepsilon}
\ee
for $e=(0,1)$\ or $(1,0)$.\\
For $\varepsilon>0$, define the function
\be
v_j(z,\varepsilon):=\iint_{D_r} \pa^j L\eta_{\varepsilon} \cdot f(\zeta) d\sigma_{\zeta}. \nn
\ee
We obtain $v_{\theta_{n-1}}(z, \varepsilon) \in C^{n-1}(D_r)$ for a fixed $\varepsilon$ by induction.
\par From (\ref{logn}) we know that, $\zeta=z$ is a singularity of $\log|z-\zeta|$ when $|j|\geq 3$. To overwhelm the blow-up behavior near the singularity, we need the Taylor expansion (\ref{Taylor expansion}). To prevent a singularity from appearing on the boundary $\pa \Omega$, we have to enlarge the domain $D_r$ of the integral (\ref{new higher}) into a larger domain $D_R$ where the divergence theorem holds. Thus for sufficiently small $\varepsilon$,
\be
&&\pa^{e_n}v_j(z, \varepsilon)=\iint_{D_r}\pa^{e_n}(\pa^{\theta_{n-1}}L\eta_{\varepsilon}) \cdot f(\zeta)d\sigma_{\zeta} \nn \\
&=&\iint_{D_R}\pa^{e_n}(\pa^{\theta_{n-1}}L\eta_{\varepsilon}) \cdot \left(f(\zeta)-P_{n-2}[f]\right)d\sigma_{\zeta} +\iint_{D_R}\pa^{e_n}(\pa^{\theta_{n-1}}L\eta_{\varepsilon}) \cdot P_{n-2}[f]d\sigma_{\zeta}. \nn
\ee
Combining the skew symmetry (\ref{skew symmety}), Green's identity (\ref{Green's identity}) and Theorem \ref{recurrent formula}, for sufficiently small $\varepsilon$, we have
\begin{eqnarray*}
&&\iint_{D_R}\pa^{e_n}(\pa^{\theta_{n-1}}L\eta_{\varepsilon}) \cdot P_{n-2}[f]d\sigma_{\zeta} \nn \\
&=&-\iint_{D_R}\frac{\pa^{e_n}}{\pa \zeta}\pa^{e_n}(\pa^{\theta_{n-1}}L\eta_{\varepsilon}) \cdot P_{n-2}[f]d\sigma_{\zeta} \nn \\
&=&-\int_{\pa D_R} \pa^{\theta_{n-1}}L \cdot P_{n-2}[f]\langle N(\zeta), e_n \rangle |d\zeta|
+\iint_{D_R} \pa^{\theta_{n-1}}L\eta_{\varepsilon} \cdot P_{n-3}[\pa^{e_n}f]d\sigma_{\zeta}\\
&=& \ldots \\
&=&-\int_{\pa D_R} \pa^{\theta_{n-1}}L \cdot P_{n-2}[f]\langle N(\zeta), e_n \rangle |d\zeta|- \ldots
-\int_{\pa D_R} \pa^{\theta_{2}}L \cdot P_{1}[\pa^{\phi_2}f]\langle N(\zeta), e_3 \rangle |d\zeta|\\
&& +\iint_{D_R} \pa^{\theta_{2}}L\eta_{\varepsilon} \cdot P_{0}[\pa^{\phi_1}f]d\sigma_{\zeta} \\
&=&-\int_{\pa D_R} \pa^{\theta_{n-1}}L \cdot P_{n-2}[f]\langle N(\zeta), e_n \rangle |d\zeta|- \ldots
-\int_{\pa D_R} \pa^{\theta_{2}}L \cdot P_{1}[\pa^{\phi_2}f]\langle N(\zeta), e_3 \rangle |d\zeta|\\
&&-\int_{\pa D_R} \pa^{\theta_{1}}L \cdot P_{0}[\pa^{\phi_1}f]\langle N(\zeta), e_2 \rangle |d\zeta|\\
&=& - \sum^{n-1}_{\tau=1} \int_{\pa D_R} \pa^{\theta_{\tau}}L \cdot P_{\tau-1}[\pa^{\phi_{\tau}}f]\langle N(\zeta), e_{\tau+1} \rangle |d\zeta|.
\end{eqnarray*}
Therefore
\be \label{patial vz}
\pa^{e_n}v_j(z, \varepsilon)&=&\iint_{D_R}\pa^{e_n}(\pa^{\theta_{n-1}}L\eta_{\varepsilon}) \cdot \left(f(\zeta)-P_{n-2}[f]\right)d\sigma_{\zeta} \nn \\
&& -\sum^{n-1}_{\tau=1} \int_{\pa D_R} \pa^{\theta_{\tau}}L \cdot P_{\tau-1}[\pa^{\phi_{\tau}}f]\langle N(\zeta), e_{\tau+1} \rangle |d\zeta|.
\ee
Now we compare (\ref{ujz}) and (\ref{patial vz}). By the local H\"{o}lder continuity of $f$, Theorem \ref{Taylor theorem} and the estimate (\ref{error term}), there exist constants $M_1$ and $M_2$ such that
\begin{eqnarray*}
&&\left|u_j(z)-\pa^{e_n}v_j(z,\varepsilon)\right| \\
&=&\left|\iint_{|\zeta-z|\leq 2\varepsilon}\left(\pa^jL-\pa^jL\eta_{\varepsilon}\right)R_{n-1}(z,\zeta)d\sigma_{\zeta}\right| \\
&\leq & M_1 \iint_{|\zeta-z|\leq 2\varepsilon} \left( \frac{n!}{|\zeta-z|^n}+\frac{2(n-1)!}{\varepsilon|\zeta-z|^{n-1}}\right)|z-\zeta|^{\nu+n-2} d\sigma_{\zeta}\\
&=& M_1 \iint_{|\zeta-z|\leq 2\varepsilon} \left( \frac{n!}{|\zeta-z|^2}+\frac{2(n-1)!}{\varepsilon|\zeta-z|}\right)|z-\zeta|^{\nu} d\sigma_{\zeta}\\
&\leq & M_2 \cdot (2\varepsilon)^{\nu}
\end{eqnarray*}
The last inequality comes from Lemma 4.2 in \cite{PDE}. Hence $\pa^{e_n}v_j(z,\zeta)$ converges to $u_j(z)$ uniformly on any compact subset of $D_r$ as $\varepsilon\rightarrow 0$. It is easy to see $v_j(z,\varepsilon)$ converges uniformly to $\pa^{\theta_{n-1}}\omega$ in the disk $D_r$, then $\omega \in C^n(D_r)$ and $u_j(z)=\pa^j \omega(z)$. The proof is complete. \hfill $\Box$
\par We list two results on a class of conformal metrics with negative curvatures as an application of potential theory. No proof is involved here. For more details, see \cite{Rothbehaviour, Rothhyper}, also \cite{Zhang1}.
\begin{theorem} $\mathrm{[6]}$ \label{original}
\textsl{ Let $\kappa:\mathbb{D}\rightarrow\mathbb{R}$ be a locally H\"{o}lder continuous function with $\kappa(0)<0$. If $u:\mathbb{D}^*\rightarrow \mathbb{R}$ is a $C^2$-solution to $\Delta u=-\kz e^{2u}$ in $\mathbb{D}^*$, then $u$ has the order $\alpha \in (-\infty, 1]$ and
\begin{align}
&u(z)=-\alpha\log|z|+v(z),            & \textrm{if\ } \ \alpha<1,\nonumber \\
&u(z)=-\log|z|-\log\log(1/|z|)+w(z),  & \textrm{if\ } \ \alpha=1,\nonumber
\end{align}
where the remainder functions $v(z)$ and $w(z)$ are continuous in $\mathbb{D}$. Moreover, the second partial derivatives satisfy the following,
\begin{align}
&v_{zz}(z), v_{z\bar{z}}(z) \ \textrm{and}\  v_{\bar{z}\bar{z}}(z) \ \textrm{are\  continuous \ at}\ z=0 & \textrm{if\ }& \ \alpha \leq 0;\nonumber\\
& v_{zz}(z), v_{z\bar{z}}(z), v_{\bar{z}\bar{z}}(z)= \textit {O}(|z|^{-2\alpha}) &\textrm{if\ }&\ 0<\alpha<1,\nonumber\\
& w_{zz}(z), w_{\bar{z}\bar{z}}(z), w_{z\bar{z}}(z)=\textit {O}(|z|^{-2}\log^{-2}(1/|z|)) & \textrm{if\ }& \ \alpha=1, \nn
\end{align}
when $z$ tends to $z=0$.}
\end{theorem}
\begin{theorem} $\mathrm{[11]}$ \label{estimate v}
\textsl{ Let $\kappa:\mathbb{D}\rightarrow\mathbb{R}$ satisfy $\kappa(0)<0$, $\kz \in C^{n-2,\nu}(\mathbb{D}^*)$ for an integer $n \geq 3$, $0<\nu \leq 1$ and let $u:\mathbb{D}^*\rightarrow \mathbb{R}$ be a $C^{n,\nu}$-solution to $\Delta u=-\kz e^{2u}$ in $\mathbb{D}^*$. Then $u(z)$ has a singularity at the origin of the order $0< \alpha \leq 1$, and for $n_1,\ n_2 \geq 1$, $n_1+n_2 =n$, near the origin, $v(z)$, $w(z)$ as in Theorem \ref{original} satisfy
$$\partial^n v(z), \ \bar{\partial}^n  v(z), \ \bar{\partial}^{n_1}\partial^{n_2}v(z)=\textit{O}(|z|^{2-2\alpha-n}),$$
$$\bar{\partial}^n w(z), \ \partial^nw(z)=\textit {O}(|z|^{-n}\log^{-2}(1/|z|)),$$
$$\bar{\partial}^{n_1}\partial^{n_2}w(z)=\textit {O}(|z|^{-n}\log^{-3}(1/|z|)),$$
where $$\partial^n=\frac{\partial ^n}{\partial z^n},\ \bar{\partial}^n=\frac {\partial^n}{\partial\bar{z}^n}$$ for a positive natural number
$n$.}
\end{theorem}
\vspace{4mm}
\textbf{Acknowledgement.} The author would like to thank Prof. Toshiyoki Sugawa and Prof. Le Hung Son for their helpful comments, suggestion and encouragement.\\

\providecommand{\bysame}{\leavevmode\hbox to3em{\hrulefill}\thinspace}
\providecommand{\MR}{\relax\ifhmode\unskip\space\fi MR }
\providecommand{\MRhref}[2]{%
  \href{http://www.ams.org/mathscinet-getitem?mr=#1}{#2}
}
\providecommand{\href}[2]{#2}

\newpage
\end{document}